\documentclass[12pt]{amsart}

\usepackage{graphicx, graphics}

\usepackage{hyperref}

\usepackage{mathrsfs}

\usepackage{amsmath}
\usepackage{amsthm}
\usepackage{amsfonts}
\usepackage{color}
\usepackage[english]{babel}
\usepackage[margin=1.3in]{geometry}

\usepackage{enumerate}
\usepackage{subcaption}

        \definecolor{brown}{rgb}{1,0,1}

\usepackage[latin9]{inputenc}
\usepackage{amsmath}
\usepackage{amsfonts}
\usepackage{amssymb}
\usepackage{amsthm}

\numberwithin{equation}{section}

\newcommand{\R}{{\mathbb R}}

\newcommand{\Z}{{\mathbb Z}}

\renewcommand{\phi}{\varphi}

\newcommand{\el}{\mathscr{L}}

\newcommand{\eT}{\mathscr{T}}

\newcommand{\ep}{\varepsilon}

\newtheorem{theo}{{\sc Theorem}}

\newtheorem{lem}[theo]{{\sc Lemma}}

\newenvironment{rem}{\medskip\noindent{\it Remark:\/} }{\medskip}

\title[ Robin spectral rigidity of domains with a reflectional symmetry ]
{Robin spectral Rigidity of strictly convex domains with a reflectional symmetry  }

\author{Hamid Hezari }
\address{Department of Mathematics, UC Irvine, Irvine, CA 92617, USA} \email{hezari@math.uci.edu}

\begin{document}
 \begin{abstract} This is a note on a recent paper of De Simoi - Kaloshin - Wei \cite{DKW}. We show that using their results combined with wave trace invariants of Guillemin-Melrose \cite{GM2} and the heat trace invariants of Zayed \cite{Za} for the Laplacian with Robin boundary conditions, one can extend the Dirichlet/Neumann spectral rigidity results of \cite{DKW} to the case of Robin boundary conditions. We will consider the same generic subset as in \cite{DKW} of smooth strictly convex $\Z_2$-symmetric planar domains sufficiently close to a circle, however we pair them with arbitrary $\Z_2$-symmetric smooth Robin functions on the boundary and of course allow deformations of Robin functions as well. 
\end{abstract}

\maketitle

\section{Introduction}

In \cite{DKW}, it is shown that for a generic class $\mathcal C$ of smooth strictly convex $\Z_2$-symmetric planar domains sufficiently close to a circle, endowed with Dirichlet or Neumann boundary conditions, one has Laplace spectral rigidity within $\mathcal C$.  This means that given any $\Omega_0 \in \mathcal C$ and any $C^1$-deformation $\{\Omega_s\}_{s \in [0, 1]} $ of $\Omega_0$ in $\mathcal C$ with $\text{Spec}(\Delta_{\Omega_s}) = \text{Spec}(\Delta_{\Omega_0})$ for all $s \in [0, 1]$, one can find isometries $\{ \mathcal I_s\}_{s \in [0, 1]}$ of $\R^2$ such that $ \mathcal I_s( \Omega_0) = \Omega_s$.  Here $\text{Spec}(\Delta_{\Omega})$ is the spectrum of the euclidean Laplacian $\Delta= \frac{\partial^2}{\partial x^2}+\frac{\partial^2}{\partial y^2}$ with Dirichlet (or Neumann) boundary condition on $\Omega$. In this paper we are concerned with the generalization of this problem for $\text{Spec}(\Delta_{\Omega, K})$ i.e.,  the spectrum of the euclidean Laplacian with Robin boundary condition $\partial_n u = Ku$ on $\partial \Omega$, for a given function $K \in C^\infty (\partial \Omega)$, where $\partial_n$ is the inward normal differentiation. In particular, by this notation $\Delta_{\Omega, 0}$ is the Laplacian on $\Omega$  with Neumann boundary condition on $\partial \Omega$. We show that: 

\begin{theo}\label{main} Let $\delta >0$ and $\mathcal S_\delta$ be the class of smooth strictly convex $\Z_2$-symmetric\footnote{It means that there exists a reflection across a line in $\R^2$ that preserves the domain.} planar domains that are $\delta$-close \footnote{See definition  2.9 of \cite{DKW}.} to a circle. Then there exists $\delta>0$ and a generic subset $\mathcal C$ of $\mathcal S_\delta$ such that given any $\Omega_0 \in \mathcal C$  and $K_0 \in C^\infty_{\Z_2} (\partial \Omega)$, and any $C^1$-deformation $\{\Omega_s\}_{s \in [0, 1]}$ of $\Omega_0$ in $\mathcal C$ and $C^0$-deformation $\{K_s\}_{s \in [0, 1]}$ of $K_0$ in $C_{\Z_2}^\infty (\partial \Omega)$ satisfying $\text{Spec}(\Delta_{\Omega_s, K_s}) = \text{Spec}(\Delta_{\Omega_0, K_0})$ for all $s \in [0, 1]$, one can find isometries $\{ \mathcal I_s\}_{s \in [0, 1]}$ of $\R^2$ such that $ \mathcal I_s( \Omega_0) = \Omega_s$ and $K_s ( \mathcal I_s (b))= K_0(b)$ for all $b \in \partial \Omega_0$. 
\end{theo}
Here, $C^\infty_{\Z_2} (\partial \Omega)$ is the space of smooth functions on $\partial \Omega$ that are invariant under the imposed $\Z_2$-symmetry on $\Omega$. Also, in fact the generic class $\mathcal C$ consists of $\Omega \in \mathcal S_\delta$ that satisfy:
\begin{itemize} \label{12}
\item[(1)] Up to the reflection symmetry, all distinct periodic billiard orbits in $\Omega$ have distinct lengths.
\item[(2)] All (transversal) periodic billiard orbits in $\Omega$ are non-degenerate, i.e., the linearized Poincar\'e map associated to each orbit does not have $1$ as an eigenvalue.
\end{itemize}

 Using the results of \cite{PSgeneric} one sees that $\mathcal C$ is generic\footnote{Countable intersection of open dense subsets with respect to Whitney $C^\infty$ topology. See \cite{PSgeneric}.}  in $\mathcal S_\delta$. Moreover, for every $\Omega \in \mathcal C$, the spectrum of $\Delta$ with Dirichlet, Neumann, or Robin boundary conditions, determines the length spectrum $\text{LS}(\Omega)$, which is the set of lengths of periodic billiard trajectories and their iterations also including the length of the boundary and its multiples with positive integers. Such determination is shown through the so called \em Poisson relation \em proved by \cite{AnMe, PS}, which asserts that if the boundary of $\Omega$ is smooth then
$$ \text{SingSupp} \left ( \text{Tr} \; \cos{ t \sqrt{-\Delta^B_\Omega} } \right ) \subset \{0 \} \cup  \pm \text{LS}(\Omega),$$ where $\Delta^B_\Omega$ is the Euclidean Laplacian with Dirichlet, Neumann, or Robin boundary conditions. One can see (\cite{PS, PSgeneric}) that under the generic conditions (1) and (2) above, the containment in the Poisson relation is an equality, hence $\text{LS}(\Omega)$ is a spectral invariant.  On the other hand, the length spectral rigidity result of \cite{DKW} shows that if $\Omega_s \in \mathcal S_\delta$, and if $\text{LS}(\Omega_s) =\text{LS}(\Omega_0)$, then  there exist isometries $\{ \mathcal I_s\}_{s \in [0, 1]}$ of $\R^2$ such that $ \mathcal I_s( \Omega_0) = \Omega_s$. Hence Theorem \ref{main} follows from the second part of the following theorem which concerns a fixed domain. To present the statement it is convenient to  fix the axis of symmetry and also a marked point as in \cite{DKW}; we assume that each $\Omega \in \mathcal S_\delta$ is invariant under the reflection about the $x$-axis, that $\Omega \subset \{ (x, y); x \geq 0 \}$, and that $0=(0, 0) \in \partial \Omega$, which will be called the marked point. 
\begin{theo}\label{Robin} Let $\mathcal C \subset \mathcal S_\delta$  be defined as above. There exists $\delta>0$ such that 
\begin{itemize}
\item[(a)] If $\Omega \in \mathcal C$,  $K_1, K_2 \in C^\infty_{\Z_2}(\partial \Omega)$, $K_1(0)=K_2(0)$, and $\text{Spec}(\Delta_{\Omega, K_1}) = \text{Spec}(\Delta_{\Omega, K_2})$, then $K_1=K_2$.
 
\item[(b)] If $\Omega \in \mathcal C$ and if there are three functions $K_1, K_2, K_3$ in $C^\infty_{\Z_2}(\partial \Omega)$ such that $$\text{Spec}(\Delta_{\Omega, K_1})= \text{Spec}(\Delta_{\Omega, K_2})= \text{Spec}(\Delta_{\Omega, K_2}),$$ then at least two of them are identical. 
\end{itemize}
\end{theo}
One can see that if we add the assumption that $\Omega$ has two perpendicular reflectional symmetries, and $K_1$ and $K_2$ are preserved under both symmetries, then $K_1(0)=K_2(0)$, hence $K_1=K_2$ by part (a).  As a result one gets the following extension of the inverse spectral result of Guillemin-Melrose \cite{GM1} obtained on ellipses. 
\begin{theo} \label{2symmetries} Let $\mathcal S_{2, \delta}$ be the subclass of $\mathcal S_{\delta}$ consisting of domains with two reflectional symmetries whose axes are perpendicular to each other. Let $\mathcal C_2 \subset \mathcal S_{2, \delta}$ be the class of domains satisfying the generic properties (1) and (2) above.  If $\Omega \in \mathcal C_2$,  $K_1, K_2 \in C^\infty_{\Z_2 \times \Z_2}(\partial \Omega)$, and $\text{Spec}(\Delta_{\Omega, K_1}) = \text{Spec}(\Delta_{\Omega, K_2})$, then $K_1=K_2$. 
\end{theo}

To prove Theorem \ref{Robin} we will use some technical results from \cite{DKW}. To be able to do so we  will need a sufficient  number of spectral invariants which we will obtain from a Poisson summation formula of Guillemin-Melrose \cite{GM2}, and also heat trace formulas of Zayed \cite{Za} for the Robin Laplacian. In fact to our knowledge these are the only Robin spectral invariants that are explicitly given in the literature. We will review these trace invariants in the next section.

\textbf{Historical background. } There is a huge literature on inverse spectral problems. Here we shall only mention the positive results concerning smooth euclidean domains and refer the reader to the surveys \cite{Z, MeSurvey, DaHe} for further historical background on positive results, and for negative results (counterexamples) we refer to \cite{GWW} and the surveys \cite{Go} and \cite{GPS}. 

Kac \cite{Ka} proved that disks in $\R^n$ are spectrally determined among all other domains. Marvizi-Melrose \cite{MM} showed that there exists a two-parameter family of smooth strictly convex domains in $\R^2$ with the symmetries of the ellipse that are spectrally isolated in an open dense class of smooth strictly convex domains. Melrose \cite{Me} and Osgood-Phillips-Sarnak \cite{OPS} established compactness of isospectral sets of smooth planar domains. Colin de Verdi\`ere \cite{CdV} proved that real analytic planar domains with the symmetries of the ellipse are spectrally rigid (i.e. all isopsectral deformations are trivial) among themselves. Zelditch \cite{Ze} proved that generic real analytic planar domains with one reflectional symmetry are spectrally distinguishable from one another. In \cite{HeZe}, it was shown that real analytic domains in $\R^n$ with reflectional symmetries about all coordinate axes are spectrally determined among the same domains. In \cite{HeZe}, ellipses were shown to be to infinitesimally spectrally rigid among smooth domains with the symmetries of the ellipse (see \cite{ PT1, PT2, PT3} for results in the context of completely integrable tables other than ellipses). Guillemin-Melrose \cite{GM1} showed Robin spectral rigidity of ellipses when the Robin functions preserve both reflectional symmetries. To our knowledge, Theorem \ref{main} is the first inverse spectral result that allows both the boundary and the Robin function to vary. The new feature is that beside the underlying domain one can also determine an additional data (namely $K$) from the trace invariants.

\section{Trace invariants for the Robin Laplacian}

\subsection{Wave trace invariants} The following is a more precise form of a result of \cite{GM2}.

\begin{theo}[Guillemin-Melrose \cite{GM2}] \label{GM}  Let $\Omega$ be a smooth strictly convex planar domain. Let $\gamma$ be a periodic billiard trajectory of length $T$, and $\{b_j \}_{j=1}^q$ be the points of reflections of $\gamma$ on $\partial \Omega$ and $\{\phi_j\}_{j=1}^q$ in $(0, \frac \pi2]$ be the angles of reflections with respect to the tangent lines at $\{b_j \}_{j=1}^q$ . Assume no\footnote{Up to reflectional symmetries when the domain is $\Z_2$ or $\Z_2 \times \Z_2$ symmetric. } periodic orbit in $\Omega$ other than $\gamma$ has length $T$, and that the linearized Poincar\'e map $\mathcal P_\gamma$ associated to $\gamma$ does not have $1$ as an eigenvalue.  Then for any $K \in C^\infty(\partial \Omega)$ we have the following singularity expansion near $t=T$
\begin{equation} 
\small Tr \left ( \cos\big( t \sqrt{-\Delta_{\Omega, K}}\big) -\cos\big( t \sqrt{-\Delta_{\Omega, 0}}\big) \right ) \sim \text{Re} \; \left ( \log (t -T+ i 0^+ ) \sum_{k=0}^\infty c_k (t-T)^k \right ),
\end{equation}
where
\begin{equation} 
c_0= C_\gamma \sum_{j=1}^q \frac{K(b_j)}{\sin(\phi_j)},
\end{equation}
for a certain constant $C_\gamma$ that depends only on $\gamma$ and is independent of $K$. Here, $\log (t -T+ i 0^+ )$ is the distribution defined by $ \lim_{\varepsilon \to 0^+ }\log (t -T+ i \varepsilon)$.
\end{theo}

\subsection{Heat trace invariants} The following heat trace formula is the main result of \cite{Za}.
\begin{theo}[Zayed \cite{Za}] \label{Za} Let $\Omega$ be a smooth simply connected planar domain and $K \in C^\infty(\partial \Omega)$.  Let $\sigma$ be an arc-length parametrization of $\partial \Omega$ in the counter-clockwise direction and $\kappa(\sigma)$ be the curvature of $\partial \Omega$ at $\sigma$. Then as $t \to 0^+$
\begin{equation}\label{Zinvariants}
\small Tr \left ( e^{t\Delta_{\Omega, K}} - e^{t\Delta_{\Omega, 0}} \right ) = \frac{1}{2 \pi} \int_{\partial \Omega} K(\sigma) d\sigma + \frac{\sqrt{t}}{ 8 \sqrt{\pi}} \int_{\partial \Omega} \left ( K(\sigma) \kappa(\sigma)+2 K^2(\sigma) \right )d\sigma+ O(t).  
\end{equation}

\end{theo}

\section{Proofs of Theorems \ref{main}, \ref{Robin}, and \ref{2symmetries}} As we discussed in the introduction, using the length spectral rigidity result of \cite{DKW}, Theorem \ref{main} reduces to part (b) of Theorem \ref{Robin}, hence we only need to prove Theorem \ref{Robin}.  

Since $\Omega \in \mathcal C$, all periodic billiard orbits $\gamma$ satisfy the required conditions of Theorem \ref{GM}. Therefore,  if we let $K=K_1-K_2$ we have
\begin{equation}\label{K}
\sum_{j=1}^q \frac{K(b_j(\gamma))}{\sin(\phi_j(\gamma))} =0,
\end{equation}
for all periodic billiard orbits. The equation \eqref{K} is very similar to (but not the same as) the equations \eqref{variation}-\eqref{lq} below, which were studied in  \cite{DKW}.  Let us first review the ingredients we need from their article. 

\subsection*{The method of De Simoi-Kaloshin-Wei}
We recall from the introduction that we have assumed that each $\Omega \in \mathcal S_\delta$ is invariant under the reflection about the $x$-axis, $\Omega \subset \{ (x, y); x \geq 0 \}$, and $0=(0, 0) \in \partial \Omega$, which is called the marked point.  At the first step, the authors show that for any $q \geq 2$, there exists a $\Z_2$-symmetric $q$-periodic orbit of rotation number $\frac{1}{q}$ passing through the marked point and having maximal length among $q$-periodic orbits of rotation number $\frac{1}{q}$ passing through the marked point. They call such an orbit a \em marked symmetric maximal $q$-periodic orbit \em and denote its length by $\Delta_q$. Next, suppose $\{\Omega_s \}$ is a $C^1$ deformation of $\Omega_0$ in $\mathcal S_\delta$ that fixes the length spectrum i.e., $\text{LS}(\Omega_s)=\text{LS}(\Omega_0)$ for all $s \in [-1, 1]$. For simplicity we represent $\Omega_s$ using a $C^1$ family $\rho_s \in C_{\Z_2}^\infty$ so that $$\partial \Omega_s = \{ b+ \rho_s(b) n(b) \, ; \, b \in \partial \Omega_0 \}, $$  where $n(b)$ is the unit inward normal at $b$. Then by taking the variation of $\Delta_q(s)$, the length of a marked symmetric maximal $q$-periodic orbit in $\Omega_s$, they show that 
\begin{equation} \label{variation}
\ell_q( \dot \rho) =0, \qquad q\geq 2,
\end{equation}
where $\dot \rho =\frac{d}{ds}|_{s=0} (\rho_s)$, and the functional $\ell_q$ is defined by
\begin{equation} \label{lq}
\ell_q(u) =\sum_{j=1}^q u(b_j)\sin(\phi_j),  \qquad q \geq 2,
\end{equation}
where $\{b_j \}_{j=1}^q$ are the points of reflections of a marked symmetric maximal $q$-periodic orbit in $\Omega_0$ and $\{\phi_j\}_{j=1}^q  \subset (0, \frac \pi2]$ are the corresponding angles of reflections. Then the authors define the map 
$$ \begin{cases} \mathcal T: C^\infty_{\Z_2} (\partial \Omega_0) \to \ell^\infty, \\ {\mathcal T}(u) =  ( {\ell}_q(u))_{q=0}^\infty, \end{cases} $$ where $\ell_q$ for $q \geq 2$ are defined above and for $q=0$ and $q=1$ are defined by
$$ {\ell}_0(u) = \int_{\partial \Omega_0} \frac{u(\sigma)}{ \kappa(\sigma)} d\sigma, $$
$$ {\ell}_1(u) =  \mu(0)u(0).$$
We recall that $\sigma$ is the arc-length parametrization of $\partial \Omega_0$, identifying $\sigma=0$ with the marked point, and $\kappa(\sigma)$ is its curvature. Also, $\mu(0)$ is the value of the Lazutkin weight $\mu$ (to be defined later) at the marked point. Since the marked point is fixed through the deformation we know that $\ell_1( \dot \rho)=0$. In addition, by taking the variation of the length of the boundary $\partial \Omega_s$ the authors show that $\ell_0 ( \dot \rho)=0$. Hence by \eqref{variation} we have 
$$ \mathcal T (\dot \rho) =0.$$  The main result of \cite{DKW} is to show that $\mathcal T$ is injective. To do this they take advantage of the Lazutkin coordinate in which they provide a precise description of periodic orbits creeping along the boundary. The Lazutkin coordinate is defined in terms of $\sigma$ by  ($\sigma =0$ being the marked point) 
$$ x( \sigma) =  C_L \int_0^ \sigma  \kappa(\sigma)^{-2/3} d \sigma, \qquad C_L= \left (  \int_{\partial \Omega_0} \kappa(\sigma)^{-2/3} d \sigma \right ) ^{-1}. $$  
In this coordinate, the Lazutkin weight is defined by 
$$ \mu(x)=  \frac{\kappa(x)^{-1/3}}{2 C_L}. $$ 
We note that for every $\ep >0$ there exists $\delta >0$ so that for every $\Omega \in \mathcal S_\delta^r$ we can make sure that $\|\mu(x) -\pi \|_{C^0} < \ep$ and  $\| \mu^{(m)} \|_{C^0} < \ep$ for $1 \leq m \leq r$ for any fixed $r \in \mathbb N$. One can easily check that for the unit circle $\mu(x)=\pi$.  

The following lemma is a crucial ingredient in \cite{DKW}, which is also a key for our proof. 
\begin{lem}\label{MainLemma} Assume $r \geq 8$.  For any $\ep>0$ sufficiently small, there exists $\delta >0$  so that for any $\Omega \in \mathcal S_\delta$ there exist $C^{r-4}$ real valued functions $\alpha(x)$ odd and $\beta(x)$ even so that for any maximal marked symmetric $q$-periodic orbit $\gamma$: 
$$ x_q^k =\frac{k}{q}+ \frac{\alpha(k/q)}{q^2} + \varepsilon O(q^{-4}),$$
$$ \phi_q^k =\frac{\mu(x_q^k)}{q} \left (1+ \frac{\beta(k/q)}{q^2} + \varepsilon O(q^{-4}) \right ),$$
where $\| \alpha \|_{C^{r-4}} =O_r(\ep)$ and $\| \beta \|_{C^{r-4}} =O_r(\ep)$. Here $\{x^k_q \}_{k=1}^q$ are the points of reflections of $\gamma$ on $\partial \Omega$ in the Lazutkin coordinate and $\{\phi^k_q\}_{k=1}^q$ in $(0, \frac \pi2]$ are the corresponding angles of reflections. 
\end{lem}
To prove their theorem, the authors use this lemma to show that the operator $\widetilde {\mathcal T}$ defined by $\widetilde {\mathcal T} (u) = \mathcal T (\frac{u}{\mu})$ is injective, whose injectivity is clearly equivalent to the injenctivity of $\mathcal T$. 

\begin{rem} We will follow the same approach expect that in our case the factor $sin(\phi_q^k)$ will be in the denominator. Another difference is that in our case we do not know the vanishing of the operators $\ell_0$ and $\ell_1$ at $K_1-K_2$. However, as we shall see, we can obtain the vanishing of $\ell_0$ by taking the limit of $\ell_q$ as $q \to \infty$. We will overcome the the lack of $\ell_1$ by using a heat trace invariant of \cite{Za}. 
\end{rem}

\subsection{Proof of part (a) of Theorem  \ref{Robin}} Throughout this section we assume that $K_1(0)=K_2(0)$ and we denote $K=K_1-K_2$. Hence by this notation $K(0)=0$. 
First we define  the operator $ \eT (u) = ({\el}_q(u))_{q=0}^\infty$, where
$$ q \geq 2: \qquad  {\el}_q(u) = \sum_{k=0}^{q-1} u(x_q^k) \frac{\mu(x_q^k)}{q^2\sin \phi_q^k},$$
$$ \el_0(u) = \int_0^1 u(x) dx, $$
$$ \el_1(u) = u(0).$$
Then it is clear from \eqref{K} that 
$$ q\geq 2:  \qquad \el_q \left ( \frac{K}{\mu} \right )=0.$$ 
Also since by our assumption $K(0)=0$, we have $$\el_1 \left ( \frac{K}{\mu} \right ) =0.$$ 
On the other hand using Lemma \ref{MainLemma}, the mean value theorem, the approximation $\sin x =x +O(x^3)$, and the Riemann sum definition of integrals, we have
$$ \lim_{q \to \infty} \el_q (u) = \el_0( u), $$
which implies that $ \el_0 \left ( \frac{K}{\mu} \right ) =0.$ As a result we have 
$$ \eT \left ( \frac{K}{\mu} \right ) =0. $$ To conclude part (a) of Theorem \ref{Robin} we need to show that $\eT$ is injective. 

We first simplify $\mu(x^k_q)$ using the asymptotic of $x^k_q$. By the mean value theorem and using the fact that $\mu(x)$ has a uniform positive lower bound
\begin{equation} \label{mu} \mu(x_q^k) = \mu \left (\frac{k}{q}+ \frac{\alpha(k/q)}{q^2} \right  ) \left ( 1 +   \varepsilon O(q^{-4}) \right ). \end{equation} Plugging this into the equation of $\phi^k_q$ we get 
$$ \phi_q^k = \frac{1}{q} \mu \left (\frac{k}{q}+ \frac{\alpha(k/q)}{q^2} \right  )\left (1+ \frac{\beta(k/q)}{q^2} + \varepsilon O(q^{-4}) \right ). $$ Next, we take $sin$ of both sides, use the mean value theorem again and the lower bound $\sin(\frac{\mu(x)}{q}) \geq \frac{C}{q}$ to get
$$ \sin \phi_q^k =  \sin \left ( \frac{1}{q} \mu \left (\frac{k}{q}+ \frac{\alpha(k/q)}{q^2} \right  )\left (1+ \frac{\beta(k/q)}{q^2} \right ) \right ) \left ( 1+ \varepsilon O_r(q^{-4}) \right ). $$
 Combining with \eqref{mu}, we obtain
$$ \frac{\mu(x^k_q)} { q^2 \sin \phi_q^k} = \frac { \frac{1}{q} \mu \left (\frac{k}{q}+ \frac{\alpha(k/q)}{q^2} \right  )}  { q \sin \left ( \frac{1}{q} \mu \left (\frac{k}{q}+ \frac{\alpha(k/q)}{q^2} \right  ) \left (1+ \frac{\beta(k/q)}{q^2} \right ) \right )} \left ( 1+ \varepsilon O_r(q^{-4}) \right ). $$
Since $$ \left (1 + \frac{\beta(k/q)}{q^2} \right )^{-1} = 1 - \frac{\beta(k/q)}{q^2} + \varepsilon O_r(q^{-4}), $$
we can rewrite the above expression as
$$ \frac{\mu(x^k_q)} { q^2 \sin \phi_q^k} = \frac { \frac{1}{q} \mu \left (\frac{k}{q}+ \frac{\alpha(k/q)}{q^2} \right  )  \left (1 + \frac{\beta(k/q)}{q^2} \right ) }  { q \sin \left ( \frac{1}{q} \mu \left (\frac{k}{q}+ \frac{\alpha(k/q)}{q^2} \right  ) \left (1+ \frac{\beta(k/q)}{q^2} \right ) \right )} \left ( 1 - \frac{\beta(k/q)}{q^2} + \varepsilon O_r(q^{-4}) \right ). $$
Applying the mean value theorem to the principal term, using $ \| \alpha \|_{C^0}, \|\beta \|_{C^0} \leq  C \varepsilon$, and 
$$\frac{x}{\sin x}= 1 + \frac{x^2}{6} + O(|x|^4), \qquad |x| < \frac{\pi+\varepsilon}{2},$$ 
we get 
$$ \frac   {\mu(x^k_q)}  { q^2 \sin \phi_q^k}=  \frac{1}{q} \left (\frac { \frac{1}{q} \mu \left (\frac{k}{q} \right  )}{ \sin \left ( \frac{1}{q} \mu \left (\frac{k}{q} \right  ) \right )} - \frac{\beta(k/q)}{q^2} + \varepsilon O_r(q^{-4}) \right ).$$
We shall write this expression in the form 
\begin{equation}\label{Sq} \frac {\mu(x^k_q)} {q^2  \sin \phi_q^k} =  \frac{1}{q} \left (1- \frac{\beta(k/q)}{q^2} + \varepsilon O_r(q^{-4}) \right ) + \frac{1}{q} S_q(\frac{k}{q}), \end{equation}
where 
\begin{equation} \label{S} S_q(x) =\frac { \frac{1}{q} \mu \left ( x \right  )} { \sin \left ( \frac{1}{q} \mu \left (x \right  ) \right )} -1. \end{equation}
\begin{rem} In \cite{DKW}, the function $S_q$ is given by $\frac { \sin \left ( \frac{1}{q} \mu \left (x \right  ) \right )}{ \frac{1}{q} \mu \left ( x \right  )} -1$. \end{rem}

To show that the operator $ \eT (u) = ({\el}_q(u))_{q=0}^\infty$ is injective we first choose a basis for $L^2_{\Z_2}( \partial \Omega)$. A convenient basis is $\{ \cos(2 \pi  j x ) \}_{j=0}^\infty$ where $x$ is the Lazutkin parameter. We shall denote $e_j=\cos(2 \pi  j x )$. Next, by some computation as performed in Lemma 5.2. of \cite{DKW}, we get
\begin{lem}\label{F} For all $q \geq 2$ and all $ j \geq 1$,
\begin{equation} \el_q(e_j) = (1- \frac{\beta_0}{q^2}) \delta_{q|j} - \frac{ \beta_j +2 \pi i j \alpha_j }{q^2}+ \mathcal S_q(e_j) +  \mathcal R _q(e_j),
\end{equation}
where the operator $\mathcal S_q$ is defined by
$$ \mathcal S _q (e_j) = \frac{1}{q} \sum_{k=0}^{q-1} S_q\left (\frac{k}{q} \right ) e_j \left (\frac{k}{q} + \frac{\alpha(\frac{k}{q})}{q^2} \right ), $$ and $\mathcal R_q$ is a remainder operator that is given by
$$ \mathcal R _q(e_j)=  \frac{1}{q^2} \sum_{s \in \Z, s \neq 0, sq \neq j }  \left (- \beta_{sq-j} + 2\pi ij  \alpha_{sq-j} \right ) + \ep O(\frac{j^2}{q^4} ).$$
The symbol $\delta_{q | j} =1$ if $q|j$ and it is zero otherwise. 
\end{lem} 

We now analyze $\mathcal S_q(e_j)$. First let us record some properties of the function $S_q$. Since in the interval $|x| < \frac{\pi+\ep}{2}$ we have 
$$ \left |\frac{x}{\sin x} -1 \right |=  \left |\frac{x-\sin x}{\sin x} \right | \leq \frac{x^3/6}{(\frac{2 \cos \ep}{\pi+\ep}) x }= \frac{(\pi+\ep) |x|^2}{12 \cos \ep},$$ we obtain the following supnorm estimates on $S_q$:
\begin{equation} \label{supnormS} | S_q(x) | \leq \frac{ (\pi +\ep) \mu^2(x)}{12q^2 \cos \ep}   \leq \frac{ (\pi +\ep)^3}{12q^2 \cos \ep}. \end{equation} Also since $| \mu'(x) | < \ep$ we have 
\begin{equation} \label{derivativesS}| S^{(r)}_q(x) | = \ep O_r( \frac{1}{q^2} ), \quad r \geq 1. \end{equation}

Now we write 
$$ \mathcal S_q (e_j) = \frac{1}{q} \sum_{k=0}^{q-1} \cos \left ( 2 \pi j \left ( \frac{k}{q} + \frac{\alpha(k/q)}{q^2} \right )\right ) S_q \left ( \frac{k}{q} \right).$$
By the mean value theorem and \eqref{supnormS} we get
$$ \mathcal S_q (e_j) = \frac{1}{q} \sum_{k=0}^{q-1} \cos \left ( \frac{2 \pi j k}{q}\right ) S_q \left ( \frac{k}{q} \right) + \ep O(\frac{j}{q^4} ).$$
We then plug in the Fourier series
$$ S_q(x) = \sum_{p \in \mathbb Z} \sigma_p(q) e^{2 \pi i p x}, $$
 of $S_q(x)$ and obtain
$$ \mathcal S_q (e_j) = \sum_{s \in \Z} \sigma_{sq-j}(q) + \ep O(\frac{j^2}{q^4} ) = \sigma_0(q) \delta_{q|j} +  \sigma_{j}(q) + \sum_{s \in \Z, s \neq 0, sq \neq j } \sigma_{sq-j}(q) + \ep O(\frac{j^2}{q^4} ). $$ 

Therefore, Lemma \ref{F} takes the following form:

\begin{lem} \label{F2} For all $q \geq 2$ and $j \geq 1$, one has
$$ {\el}_q(e_j) = (1+ \sigma_0(q) - \frac{\beta_0}{q^2}) \delta_{q|j} + \frac{\el^{*}(e_j)}{q^2} +  \mathcal R _q(e_j),
$$
where 
$$ \el^*(e_j)  =  q^2\sigma_j(q) - \beta_j -2 \pi j \alpha_j ,$$ 
and the remainder operator is given by
\begin{equation} \label{R}  \mathcal R _q(e_j)=  \frac{1}{q^2} \sum_{s \in \Z, s \neq 0, sq \neq j } q^2 \sigma_{sq-j}(q) +  \alpha_{sq-j}(q) - 2\pi ij  \beta_{sq-j}(q) + \ep O(\frac{j^2}{q^4} ) \end{equation}
\end{lem}

However, this is not a convenient way of writing the operator $ \el^* $, since $q^2\sigma_j(q)$ hence $\el^*$,  depends on $q$. To resolve this we note that
\begin{align*} \sigma_j (q) = \int_0^1 S_q(x) e^{2 \pi i j x} dx & = \int_0^1 \left ( \frac { \frac{1}{q} \mu \left ( x \right  )} { \sin \left ( \frac{1}{q} \mu \left (x \right  ) \right )} -1 \right )  e^{2 \pi i j x} dx  \\
& =  \int_0^1 \frac{ \mu^2(x)}{6q^2}  e^{2 \pi i j x} dx +  \int_0^1 R \left ( \frac{\mu(x)}{q} \right )  e^{2 \pi i j x}dx, 
\end{align*}
where $R$ is defined by $ \frac{y}{\sin(y)} - 1 = \frac{y^2}{6} +R(y).$ Since $R(y)=O(y^4)$ and $R'(y)=O(y^3)$, by performing integration by parts once to the the second integral and the fact $| \mu'(x) | < \ep$, we get
$$ \int_0^1 R \left ( \frac{\mu(x)}{q} \right )  e^{2 \pi i j x} dx = O( \frac{\ep}{jq^4}). $$
Therefore we can absorb this term in the remainder term $\mathcal R_q(e_j)$.
The conclusion is we can write 
\begin{equation}\label{Lq} \el_q(e_j) =  \left (1+ \sigma_0(q) - \frac{\beta_0}{q^2} \right ) \delta_{q|j} + \frac{\el^{**}(e_j)}{q^2} +  \mathcal R _q(e_j), \end{equation}
where
$$\el^{**}(e_j)  =  \widetilde {\sigma}_j - \beta_j -2 \pi j \alpha_j ,$$ 
with 
$$ \widetilde {\sigma}_j  = \int_0^1 \frac{ \mu^2(x)}{6}  e^{2 \pi i j x} dx. $$

We now observe that by the properties \eqref{supnormS} and \eqref{derivativesS} of $S_q(x)$ we have
\begin{equation} \label{sigma0} |\sigma_0(q)| <   \frac{ (\pi +\ep)^3}{12q^2 \cos \ep}, \end{equation}
$$ | \sigma_p(q)  | = \ep O_r( \frac{1}{j^{r}q^2} ), \quad p \neq 0. $$
Note that the second equation follows from integration by parts.
This shows that for $p \neq 0$,  $q^2\sigma_p(q)$ behaves similarly as $\alpha_p$ and $\beta_p$.

Now assume $\eT (u)=0$.  We recall that $\eT = \{ \el_q \}_{q=0}^\infty$. Then in particular 
$$\el_0(u)= \int_0^1 u(x) dx=0.$$ Thus by \eqref{Lq} and \eqref{R} we can write 
\begin{equation}\label{decomposition} \eT(u) = \el^{**}(u)   b_*+  \eT_{*, R}(u), \end{equation}
where $(b_*)_q = \frac{1}{q^2}$ for $q \geq 2$ and $(b_*)_0=(b_*)_1=0$, and $\eT_{*, R}$ is defined on the basis $\{ e_j \}_{j=0}^\infty$ by 
$$ j \geq 1: \quad \eT_{*, R}(e_j) =  \left (1+ \sigma_0(q) - \frac{\beta_0}{q^2} \right ) \delta_{q|j}  +  \mathcal R _q(e_j), $$
$$   \eT_{*, R}(e_0)=0. $$
For $ 3 < \gamma < 4$  we denote 
$$ X_\gamma = \{ u \in L^1(\mathbb T): u(x)=u(-x), \; \hat u_0=0, \; \hat u_j =o(j^{- \gamma})\},  \qquad || u ||_{X_\gamma} = \max_{ j\geq 1} j^\gamma \hat u_j \; ,$$ 
$$ \ell_\gamma= \{b \in \ell^\infty: \; b_0=0, \; b_q = o( q^{-\gamma}) \}, \qquad || b ||_{\ell_\gamma} = \max_{ j\geq 1} j^\gamma b_j \; .$$
Then one can easily see that $\eT_{*, R}$ maps $X_\gamma$ to  $\ell_\gamma$ (in fact it is invertible as we shall see below).  
If $\eT(u)=0$, then since for us $u = \dot \rho \in C^\infty( \mathbb T) \subset X_\gamma$, by \eqref{decomposition}, 
$$ \el^{**}(u) b_* \in \text{Range} ( \eT_{*, R} ) \subset \ell_{\gamma}. $$ 
However since $\gamma >3$, this is impossible unless $ \el^{**}(u)=0$, which by \eqref{decomposition} implies that
$$ {\eT}_{*, R} (u) =0.$$ Then we show that  $\eT_{*, R}: X_\gamma \to  \ell_\gamma$ is invertible by showing that  $||\eT_{*, R} - \text{Id} ||_\gamma < 1$ where $||. ||_\gamma$ is the operator norm from $X_\gamma$ to $\ell_\gamma$. Here the identity operator $\text{Id}:X_\gamma \to \ell_\gamma$ is defined by $\text{Id}(e_j)= e'_j$ for all $j \geq 1$, where as before $\{e_j\}_{j=1}^\infty = \{ \cos (2 \pi jx) \}_{j=1}^\infty$ and  $\{e'_q\}_{q=1}^\infty$ is the standard basis for $\ell_\gamma$. For any operator $T: X_\gamma \to \ell_\gamma$ with the matrix representation $[T_{qj}]$, the operator norm is given by
$$ || T ||_\gamma = \sup_{ q \geq 1} \sum_{ j\geq 1} q^\gamma j^{-\gamma} |T_{qj}|.$$ Let $\Delta: X_\gamma \to \ell_\gamma$ be the operator with the matrix $[ \delta_{q |j} ]$.  Then 
$$ (\eT_{*, R})_{1j} = (\Delta)_{1j}=1$$
$$ q \geq 2: \qquad  ({\eT}_{*, R})_{qj} = \Delta_{qj} + \left (\sigma_0(q) - \frac{\beta_0}{q^2} \right ) \Delta_{qj} + \mathcal R_{qj}.$$
By a simple estimate (see \cite{DKW}) one gets
$$ || \Delta-I ||_\gamma< \zeta(3) -1$$
In particular $||\Delta||_\gamma < \zeta(3)$. This shows that if $\Delta'$ is the operator defined by
$$ (\Delta')_{1j} =0$$
$$ q \geq 2: \quad (\Delta')_{qj} =\left (\sigma_0(q) - \frac{\beta_0}{q^2} \right ) \Delta_{qj},$$
then by the estimate \eqref{sigma0} on $\sigma_0(q)$ 
$$|| \Delta'||_\gamma \leq \left ( \frac{(\pi +\ep)^3}{48 \cos \ep} + \frac{C\ep}{4} \right ) \zeta(3).$$
On the other hand by the computations of \cite{DKW} we know that $|| \mathcal R ||_\gamma \leq C \ep$. Combining these estimates together we obtain
\begin{align*} ||\eT_{*, R} - \text{Id} ||_\gamma & = ||\Delta+ \Delta' + \mathcal R - \text{Id} ||_\gamma  \\ & \leq  ||\Delta - \text{Id} ||_\gamma + || \Delta'||_\gamma + ||\mathcal R||_\gamma \\ 
& \leq \zeta(3)-1 +  \left ( \frac{(\pi +\ep)^3}{48 \cos \ep} + \frac{C\ep}{4} \right ) \zeta(3) + C \ep.  \end{align*}
At $\ep=0$, the last expression is less than $\frac{979}{1000}$, hence by choosing $\ep>0$ sufficiently small we can guarantee that it is less than one. 

\subsection{Proof of part (b) of Theorem \ref{Robin}} 
If two of $K_1(0), K_2(0),K_3(0)$ agree then by part (a) of Theorem \ref{Robin} we are done. So assume  $K_1(0), K_2(0),K_3(0)$ are three distinct numbers and define the function $f \in C^\infty_{\Z_2}(\partial \Omega)$ by 
$$f(b) = \frac {K_2(b) - K_3(b)} {K_2(0)- K_3(0)} \; , \quad b \in \partial \Omega .$$
Then we consider the functions
$$ K_{12} = (K_1 -K_2) - (K_1(0) - K_2(0)) f \;,$$
$$ K_{13} = (K_1 -K_3) - (K_1(0) - K_3(0)) f \; .$$
Obviously $K_{12}, K_{13} \in C^\infty_{\Z_2}(\partial \Omega)$ and $K_{12}(0)=K_{13}(0)=0$ because $f(0)=1$. Since  $\text{Spec}(\Delta_{\Omega, K_1})$ $= \text{Spec}(\Delta_{\Omega, K_2})= \text{Spec}(\Delta_{\Omega, K_2})$ by the notations and the discussion at the beginning of the proof of part (a) we have
$$ \eT\left ( \frac{K_{12}}{\mu} \right )=\eT \left ( \frac{K_{13}}{\mu} \right )=0. $$ However, we showed in the previous section that the operator $\eT$ is injective. Thus $K_{12}=K_{13}=0$, or equivalently 
\begin{equation}  \label{K12} K_1 -K_2 = (K_1(0) - K_2(0)) f , \end{equation}
\begin{equation}  \label{K13}  K_1 -K_3 = (K_1(0) - K_3(0)) f . \end{equation} 
On the other hand by the heat trace formula \eqref{Zinvariants}, we know that
$$\int_{\partial \Omega} K_1 \kappa +2 K_1^2 = \int_{\partial \Omega} K_2 \kappa +2 K_2^2 =\int_{\partial \Omega} K_3 \kappa +2 K_3^2. $$ 
These imply that
$$ \int_{\partial \Omega}  (K_1 -K_2) \big ( \kappa + 2 (K_1+K_2) \big ) =0, $$
$$ \int_{\partial \Omega}  (K_1 -K_3)\big  ( \kappa + 2 (K_1+K_3) \big ) =0. $$
By plugging \eqref{K12} and \eqref{K13} into these identities and dividing by $K_1(0) - K_2(0)$ and $K_1(0) - K_3(0)$ respectively we get
$$ \int_{\partial \Omega}   \big ( \kappa + 2 (K_1+K_2) \big ) f =0, $$
$$ \int_{\partial \Omega}   \big  ( \kappa + 2 (K_1+K_3) \big ) f  =0. $$
By subtracting these two equations we obtain
$$ \int_{\partial \Omega}   \big ( K_2 - K_3 ) f =0.$$
Recalling the definition of $f$ we get $\int_{\partial \Omega}  f^2 =0.$ Since $f$ is continuous and real-valued this implies that $f=0$. However this contradicts $f(0)=1$.  

\subsection{Proof of Theorem \ref{2symmetries}} Let $0'$ be the other point on the axis of symmetry other than the marked point $0$, and consider the $2$-orbit bouncing between $0$ and $0'$. Then by \eqref{K} we have $ K_1(0) + K_1(0')= K_2(0) + K_2(0')$. On the other hand since $K_1$ and $K_2$ are $\Z_2 \times \Z_2$-symmetric we have $K_1(0)=K_1(0')$ and  $K_2(0)=K_2(0')$ . These together show that $K_1(0)=K_2(0)$, and hence by part (a) of Theorem \ref{Robin} we must have $K_1=K_2$.

\end{document}